\documentclass[11pt,a4paper]{amsart}

\usepackage{amsmath, amsthm, amssymb}

\newcommand{\bC}{\mathbb{C}}
\newcommand{\bP}{\mathbb{P}}

\newcommand{\bN}{\mathbb{N}}

\newcommand{\FS}{\operatorname{FS}}

\newcommand{\supp}{\operatorname{supp}}
\newcommand{\Id}{\mathrm{Id}}
\newcommand{\Aut}{\operatorname{Aut}}
\newcommand{\bR}{\mathbb{R}}
\newcommand{\bQ}{\mathbb{Q}}
\newcommand{\bZ}{\mathbb{Z}}
\newcommand{\bT}{\mathbb{T}}

\newcommand{\rd}{\mathrm{d}}
\newcommand{\loc}{\operatorname{loc}}
  
\numberwithin{equation}{section}
\theoremstyle{plain}
\newtheorem{theorem}{Theorem}[section]

\newtheorem{proposition}[theorem]{Proposition}

\newtheorem{mainth}{Theorem}

\theoremstyle{definition}

\newtheorem*{acknowledgement}{Acknowledgement}

\theoremstyle{remark}
\newtheorem{remark}[theorem]{Remark}

\begin{document}  
  
\title[An a priori bound of endomorphisms of $\bC\bP^k$]{An a priori bound of endomorphisms of $\bC\bP^k$ and a remark on the Makienko conjecture in dimension one}

\author[Y\^usuke Okuyama]{Y\^usuke Okuyama}
\address{
Division of Mathematics,
Kyoto Institute of Technology,
Sakyo-ku, Kyoto 606-8585 Japan}
\email{okuyama@kit.ac.jp}

\date{\today}

\begin{abstract}
Let $f$ be an endomorphism of $\bC\bP^k$ of degree $>1$, and assume that
for any cyclic Fatou component $W$ of $f$ having a period $p\in\bN$,
the equilibrium measure $\mu_f$ has a positive charge on the boundary of $W$ 
if and only if $f^{-p}(W)=W$. Then
we obtain a locally uniform a priori bound of the dynamics of $f$,
which in particular yields a Diophantine-type estimate of the dynamics of $f$ 
on its domaines singuliers. We also point out that
in the case of $k=1$, the statement of our assumption is related to
both the impossibility for the Julia set of $f$ to be the boundary of lakes of Wada 
and the so called Makienko conjecture 
on the non-emptiness of the residual Julia set of $f$.
\end{abstract}


\maketitle

\section{Introduction}\label{sec:intro}

Let $f$ be a holomorphic endomorphism of $\bP^k=\bC\bP^k$,
$k\in\bN$, of algebraic degree
$d>1$ (so of topological degree $d^k$). 
By the seminal Forn{\ae}ss--Sibony \cite{FS94nato}
(see also the survey \cite{DSsurvey}),
the weak limit 
\begin{gather*}
 \mu_f:=\lim_{n\to\infty}\frac{(f^*\omega_{\FS})^{\wedge k}}{d^{kn}}
\end{gather*}
exists on $\bP^k$, so that $\mu_{f^n}=\mu_f$ on $\bP^k$ for every $n\in\bN$,
where $\omega_{\FS}$ is the Fubini-Study K\"ahler form on $\bP^k$ 
normalized as $\omega_{\FS}^{\wedge k}(\bP^k)=1$. 
Let us equip $\bP^k$ with the normalized
{\itshape chordal} metric $[x,y]_{\bP^k}(\le 1)$
 (see \eqref{eq:chordal} for the definition),
which is comparable to the metric
on $\bP^k$ induced by $\omega_{\FS}$.
The probability measure $\mu_f$,
which is called the {\itshape equilibrium measure} of $f$, 
has no masses on pluripolar subsets in $\bP^k$ and  
satisfies the $f$-balanced property $f^*\mu_f=d^k\cdot\mu_f$ on $\bP^k$.
Letting $\delta_x$ be the Dirac measure on $\bP^k$ at each $x\in\bP^k$,
for every $x\in\bP^k$ except for a pluripolar subset in $\bP^k$,
the weak convergence
\begin{gather}
 \lim_{n\to\infty}\frac{(f^n)^*\delta_x}{d^{kn}}=\mu_f
\label{eq:equidist} 
\end{gather}
holds on $\bP^k$, so $\mu_f$ is mixing (so ergodic) under $f$
(see \cite{BD01,DS08,DS10} for more details).
Let us also recall a more classical Fatou-Julia strategy.
The (first) Julia set $J(f)$ of $f$ is by definition
all points in $\bP^k$ at each of which the family $\{f^n:n\in\bN\}$ is not normal,
and is a compact subset in $\bP^k$.
The (first) Fatou set $F(f)$ of $f$ is defined by $\bP^k\setminus J(f)$,
and each component of $F(f)$ is called a Fatou component of $f$.
Both $J(f)$ and $F(f)$ are totally invariant under $f$, 
$f$ maps each Fatou component of $f$ properly to a Fatou component of $f$,
and the preimage of a Fatou component of $f$ under $f$ is the union of
(at most $d^k$) Fatou components of $f$. A Fatou component $W$ of $f$ is
said to be cyclic under $f$ if $f^p(W)=W$ for some $p\in\bN$; then
such $p$ is called {\itshape a} period of $W$ (under $f$),
and the minimal such $p$ is called the {\em exact period} of $W$ (under $f$). 
In general, we have
$\supp\mu_f\subset J(f)$,
and we have $\supp\mu_f=J(f)$ if $k=1$ (see Section \ref{sec:lift} for more details).

Our principal result is the following locally uniform {\em a priori} bound
of the dynamics of $f$. 

\begin{mainth}[an {\itshape a priori} bound]\label{th:nonlinvanish}
Let $f$ be a holomorphic endomorphism of $\bP^k$ of algebraic degree $d>1$,
and suppose that {\em (*) for any cyclic Fatou component $W$ of $f$ having 
a period $p\in\bN$, $\mu_f(\partial W)>0$ if and only if $f^{-p}(W)=W$.}
Then for every holomorphic endomorphism $g$ of $\bP^k$ of degree $>0$
and every non-empty open subset $D$ in $\bP^k$,
\begin{gather}
\lim_{n\to\infty}\frac{\sup_{y\in D}\log[f^n(y),g(y)]_{\bP^k}}{d^n+\deg g}=0.\label{eq:negative}
\end{gather}
\end{mainth}

We conclude this section with a few remarks on the conclusion \eqref{eq:negative}
and the assumption (*) in Theorem \ref{th:nonlinvanish}
since both are seemingly technical.
\begin{remark}\label{sec:rotation}
 Following Fatou \cite[Sec.\ 28]{Fatou1920deux},
 a cyclic Fatou component $W$ of $f$ having the exact period, say, 
 $p\in\bN$ is called
 a {\itshape domaine singulier} (a singular domain, or a rotation domain)
 of $f$ if $f^p:W\to W$ is injective ({\itshape singulier} in French)
so biholomorphic.
The conclusion \eqref{eq:negative} in Theorem \ref{th:nonlinvanish} for $g=\Id_{\bP^k}$ is equivalent to that {\itshape for every {\em domaine singulier} $W$ of $f$ having the exact period $p\in\bN$ and every non-empty open subset $D\Subset W$,}
 \begin{gather}
 \lim_{n\to\infty}\frac{\log\sup_{y\in D}[f^{pn}(y),y]_{\bP^k}}{d^{pn}+1}=0.\label{eq:rotation}
 \end{gather}
 This {\itshape Diophantine-type} estimate \eqref{eq:rotation} of $f^p$ on 
 a {\em domaine singulier} $W$
 has been known in \cite[Theorem 3]{OkuNonlinear}
 without assuming (*) but
 under the assumption (**) that $W$ is of maximal type 
 in that, setting $q:=\min\{j\in p\bN:f^j|W\in G_0\}$,
 where $G_0$ is the component of the closed subgroup
 generated by $f^p|W$ in the biholomorphic automorphisms group $\Aut(W)$
 of $W$ containing $\Id_W$, 
 there is a Lie groups isomorphism $G_0\to\mathbb{T}^k$
 that maps $f^q|W$ to $(e^{2i\pi\alpha_1},\ldots,e^{2i\pi\alpha_k})\in\mathbb{T}^k$
 for some $\alpha_1,\ldots,\alpha_k\in(\bR\setminus\bQ)/\bZ$; in general,
 $G_0$ at least contains a real $1$-dimensional torus $\bT^1$
 (Ueda \cite{Ueda98}). 
\end{remark}

\begin{remark}\label{th:equivalent}
In the statement in the assumption (*), one implication is obvious;
if $f^{-p}(W)=W$, then $\supp\mu_f\subset\partial W$
by (the property of pluripolar subsets and) the above description of $\mu_f$, 
so $\mu_f(\partial W)=1>0$.

In fact, it is always the case that 
{\em for every cyclic Fatou component $W$ of $f$ having 
a period $p\in\bN$, we have $\mu_f(\partial W)\in\{0,1\}$ and, 
if $\mu_f(\partial W)=0$, then we have
$\mu_f(\partial U)=0$ for every component $U$ of 
$\bigcup_{n\in\bN\cup\{0\}}f^{-n}(W)$}; 
for, by the ergodicity of $\mu_f$ under $f$ and $f^p(W)=W$, 
we have not only
$\mu_f(\bigcup_{n\in\bN\cup\{0\}}f^{-pn}(\partial W))\in\{0,1\}$ but, 
for every $n\in\bN$, also $f^{-p(n-1)}(\partial W)\subset f^{-pn}(\partial W)$. 
Then by the $f$-invariance of $\mu_f$, we compute as
$\mu_f(\partial W)=\mu_f(f^{-pn}(\partial W))
=\mu_f(f^{-p(n-1)}(\partial W))
+\mu_f(f^{-pn}(\partial W)\setminus f^{-p(n-1)}(\partial W))
=\mu_f(\partial W)+\mu(f^{-pn}(\partial W)\setminus f^{-p(n-1)}(\partial W))$,
so $\mu_f(f^{-pn}(\partial W)\setminus f^{-p(n-1)}(\partial W))=0$, 
and then $\mu_f(\partial W) 
=\mu_f(\bigcup_{n\in\bN\cup\{0\}}f^{-pn}(\partial W))\in\{0,1\}$.
Finally, if in addition $\mu_f(\partial W)=0$, then 
for every component $U$ of $\bigcup_{n\in\bN\cup\{0\}}f^{-n}(W)$, 
by the $f$-invariance of $\mu_f$
and $f^p(W)=W$, we compute as
$0=\mu_f(\partial W)=\mu_f(\partial(f^{-(p-j)}(W)))\ge
\mu_f(\partial(f^j(W)))=\mu_f(\bigcup_{n\in\bN\cup\{0\}}f^{-pn}(\partial(f^j(W))))
\ge\mu_f(\partial U)\ge 0$ for some $j\in\{0,\ldots,p-1\}$,
so $\mu_f(\partial U)=0$.
\end{remark}

\begin{remark}
In the dimension $k=1$ case,
we have $\supp\mu_f=J(f)$, and by the Riemann-Hurwitz formula, 
for any cyclic Fatou component $W$ of $f$ having a period $p\in\bN$,
if $f^{-p}(W)=W$, then $p\in\{1,2\}$.
Moreover, by Sullivan's no wandering domain theorem
\cite{Sullivan85}, for any Fatou component $U$ of $f$,
there is $N\in\bN\cup\{0\}$ such that $W:=f^N(U)$ is a cyclic Fatou 
component of $f$; then $\mu_f(\partial W)\in\{0,1\}$ 
and, if in addition $\mu_f(\partial W)=0$, then
$\mu_f(\partial U)=0$ (seen in Remark \ref{th:equivalent}).
By the Riemann-Hurwitz formula again, there are at most two 
Fatou components $U$ of $f$ such that $f^{-2}(U)=U$.

It seems open whether the statement in the assumption (*) is the case 
for any $f$, even in the case $k=1$. We point out the following.
\begin{proposition}
 Let $f\in\bC(z)$ be of degree $d>1$. Then the following statements
 are equivalent$;$
 \begin{enumerate}
 \item for any Fatou component $U$ of $f$, 
       $\mu_f(\partial U)>0$ if and only if $f^{-2}(U)=U$.
       \label{measure}
 \item for any Fatou component $U$ of $f$,
 $\partial U=J(f)$ if and only if $f^{-2}(U)=U$.
       \label{makienko}
 \end{enumerate}
Moreover, if the above statements are the cases, then 
so are the followings$;$ 
\begin{enumerate}
\setcounter{enumi}{2}
 \item  the Julia set $J(f)$ of $f$ cannot be the boundary of lakes of Wada.
Here a compact subset $C$ in $\bP^1$ is called the {\em boundary of
lakes of Wada} if there are at least three distinct components of
$\bP^1\setminus C$, the boundary of each of which equals the whole $C$ 
$($introduced in \cite{Yoneyama17}$)$. 
 \label{wada}
 \item $($the so called {\em Makienko conjecture} $(1990))$
the residual Julia set
 $J_0(f):=J(f)\setminus(\bigcup_{U:\text{a Fatou component of }f}\partial U)$
       of $f$
 is empty if and only if there is a Fatou component $U$ of $f$ 
       such that $f^{-2}(U)=U$.
       \label{residual}
\end{enumerate}
\end{proposition} 

Indeed, the observations in the first paragraph yield
\eqref{measure}$\Leftrightarrow$\eqref{makienko},
\eqref{makienko}$\Rightarrow$\eqref{wada}, and 
\eqref{measure}$\Rightarrow$\eqref{residual}.
For further studies, see e.g.\ 
 \cite{Qiao95,Morosawa97,SunYang03,CMMR09,HT18}.
\end{remark}

\section{Background}\label{sec:lift}
We say a function $H:\bC^{k+1}\to\bR\cup\{-\infty\}$
is log-homogeneous (of order $1$) if
for every $c\in\bC^*$, $H(cZ)=H(Z)+\log|c|$ on $\bC^{k+1}$. 
Let $\|\cdot\|$ be the Euclidean norm on $\bC^{k+1}$.
The function $\log\|\cdot\|:\bC^{k+1}\to\bR\cup\{-\infty\}$ is
continuous, plurisubharmonic, and log-homogeneous
(of order $1$). 
The complex Laplacian $\rd\rd^c$ is normalized as usual,
so in particular that
$\pi^*\omega_{\FS}=\rd\rd^c\log\|\cdot\|$ on $\bC^{k+1}\setminus\{0\}$,
where the origin $(0,\ldots,0)$ of the $\bC$-linear space
$\bC^{k+1}$ is denoted by $0$. More generally,
a function $H:\bC^{k+1}\to\bR\cup\{-\infty\}$
which is plurisubharmonic and log-homogeneous (of order $1$) is also called an
$\omega_{\FS}$-plurisubharmonic function {\itshape on} $\bP^k$;
indeed, the function $H-\log\|\cdot\|$ descends to $\bP^k$ so that
$\pi^*\bigl(\rd\rd^c(H-\log\|\cdot\|)+\omega_{\FS}\bigr)=\rd\rd^c H$ on $\bC^{k+1}\setminus\{0\}$, $\rd\rd^c(H-\log\|\cdot\|)+\omega_{\FS}\ge 0$ on $\bP^k$,
and $\bigl(\rd\rd^c(H-\log\|\cdot\|)+\omega_{\FS}\bigr)\wedge\omega_{\FS}^{k-1}$
is a probability measure on $\bP^k$.

Set $\ell(k):=\binom{k+1}{2}\in\bN$ so that 
$\bC^{k+1}\wedge\bC^{k+1}\cong\bC^{\ell(k)}$ (cf.\ \cite[\S 8.1]{Kobayashi98})
and, for simplicity,
also denote by $\|\cdot\|$ the Euclidean norm on $\bC^{\ell(k)}$.
Let $\pi:\bC^{k+1}\setminus\{0\}\to\bP^k$
be the canonical projection.
The chordal metric on $\bP^k$ is
\begin{gather}
 [x,y]_{\bP^k}:=\|Z\wedge W\|(\|Z\|\cdot\|W\|)\le 1,\quad x,y\in\bP^k,\label{eq:chordal}
\end{gather}
where $Z\in\pi^{-1}(x),W\in\pi^{-1}(y)$.

Let $f$ be a holomorphic endomorphism of $\bP^k=\bC\bP^k$ 
of algebraic degree $d>1$. A lift of $f$ is 
an ordered $(k+1)$-tuple $F=(F_0,F_1,\ldots,F_k)\in(\bC[z_0,z_1,\ldots,z_k]_d)^{k+1}$ 
of homogeneous polynomials of degrees $d$ 
in the indeterminants $z_0,\ldots,z_{k+1}$,
which is unique up to multiplication in $\bC^*$, such that
$\pi\circ F=f\circ\pi$ 
on $\bC^{k+1}\setminus\{0\}$
(and that $F^{-1}(0)=\{0\}$). 
Fix a lift $F$ of $f$. Then the uniform limit
\begin{gather}
 G^F:=\lim_{n\to\infty}\frac{\log\|F^n\|}{d^n}:\bC^{k+1}\setminus\{0\}\to\bR\label{eq:escape}
\end{gather}
exists (\cite{HP94}) and, setting $G^F(0):=-\infty$, 
the function $G^F:\bC^{k+1}\to\bR\cup\{-\infty\}$ is
continuous, plurisubharmonic,
and log-homogeneous (of order $1$).
By Ueda \cite[Theorem 2.2]{Ueda94},
the region of pluriharmonicity of $G^F$ 
coincides not only with $\pi^{-1}(F(f))$ 
but also with $\pi^{-1}(\tilde{F}(f))$, where $\tilde{F}(f)$
is the set of all points in $\bP^k$ at each of which 
the family $\{f^{n_j}:\bP^k\to(\bP^k,[x,y]_{\bP^k}):j\in\bN\}$ is normal
for {\em some} sequence $(n_j)$ in $\bN$ tending to $\infty$ as $j\to\infty$. 

\section{Proof of Theorem \ref{th:nonlinvanish}}
\label{sec:vanish}
Let $f$ be a holomorphic endomorphism of $\bP^k=\bC\bP^k$ 
of algebraic degree $d>1$, and assume (*) that
for any cyclic Fatou component $W$ of $f$ having a period $p\in\bN$,
$\mu_f(\partial W)>0$ if and only if $f^{-p}(W)=W$.

Suppose to the contrary that 
there are a holomorphic endomorphism $g$ of $\bP^k$ of degree $>0$
and a domain $D$ in $\bP^k$ 
such that \eqref{eq:negative} does not hold, so that
there is a sequence $(n_j)$ in $\bN$ tending to $\infty$ as $j\to\infty$ such that
\begin{gather}
\lim_{j\to\infty}\frac{\sup_{y\in D}\log[f^{n_j}(y),g(y)]_{\bP^k}}{d^{n_j}+\deg g}<0.\label{eq:quick}
\end{gather}
Then $D\subset F(f)$, and let $U$ be the Fatou component of $f$ containing $D$.
Since $\deg g>0$, 
we have $\lim_{j\to\infty}f^{n_{j+1}-n_j}=\Id_{g(D)}$ locally uniformly
on $g(D)$, and then there is $N\in\bN\cup\{0\}$ such that $V:=f^{n_N}(U)(=g(U))$
is a cyclic Fatou component of $f$ having the exact period, say, $p\in\bN$
and satisfies $\deg(f^p:V\to V)=1$. Then $f^{-p}(V)\neq V$ since $($the topological
degree of $f)=d^k>1$, so that under the assumption (*), we have
\begin{gather}
 \mu_f(\partial U)=0\label{eq:preimage} 
\end{gather}
(seen in Remark \ref{th:equivalent}).

Also fix a lift $G$ of $g$. By \eqref{eq:chordal} and \eqref{eq:escape},
the family $\bigl\{(\log|F^{n_j}\wedge G|)/d^{n_j}:j\in\bN\bigr\}$
is locally uniformly bounded from above on $\bC^{k+1}$.
Moreover, by \eqref{eq:chordal}, \eqref{eq:escape}, and $J(f)\neq\emptyset$,
we have $\limsup_{j\to\infty}\sup_{\{G^F=0\}}
(\log|F^{n_j}\wedge G|)/d^{n_j}\ge 0>-\infty$.
Hence by \eqref{eq:chordal}, \eqref{eq:escape}, and a version of
Hartogs lemma for a sequence of plurisubharmonic functions 
(see \cite[Theorem 4.1.9(a)]{Hormander83} or \cite[Theorem 1.1.1]{Azarin80}),
taking a subsequence of $(n_j)$ if necessary, 
there is a plurisubharmonic function $\phi$ on $\bC^{k+1}$ such that
\begin{gather*}
 \phi:=\lim_{j\to\infty}\frac{\log\|F^{n_j}\wedge G\|}{d^{n_j}+\deg g}(\le G^F)
\quad\text{in }L^1_{\loc}(\bC^{k+1},m_{2k+2}),
\end{gather*}
where $m_{2k+2}$ denotes the $(2k+2)$-dimensional
Lebesgue measure on $\bC^{k+1}\cong\bR^{2k+2}$.
The function $\phi$ is also log-homogeneous
(of order $1$).

By the log-homogeneity of both $G^F$ and $\phi$,
$\phi-G^F$ on $\bC^{k+1}\setminus\{0\}$ descends to 
a function $\bP^k\to\bR_{\le 0}\cup\{-\infty\}$, 
which is not only upper semicontinuous on $\bP^k$
but also plurisubharmonic on $U$, and
by \eqref{eq:chordal}, \eqref{eq:escape}, and \eqref{eq:quick}, 
the open subset $\{\phi-G^F<0\}$ in $\bP^k$ contains $D(\subset U)$
except for an at most $m_{2k+2}$-null subset in $D$. Hence
$U\cap\{\phi-G^F<0\}\neq\emptyset$.
By \eqref{eq:chordal}, \eqref{eq:escape},
the upper semicontinuity of $\phi-G^F$,
and a version of
Hartogs lemma for a sequence of plurisubharmonic functions 
(see \cite[Theorem 4.1.9(b)]{Hormander83}),
we also have $\phi-G^F\equiv 0$ on $J(f)$, 
so in particular $\phi=G^F$ on $\pi^{-1}(\partial U)$.

Let us define the {\itshape locally bounded} function 
\begin{gather*}
 \psi:=\begin{cases}
	\max\{\phi,G^F-1\} & \text{on }\pi^{-1}(U)\\
	G^F & \text{on }\bC^{k+1}\setminus \pi^{-1}(U), 
       \end{cases}
\end{gather*}
which is still plurisubharmonic on $\bC^{k+1}$;
for, it is upper semicontinuous on $\bC^{k+1}$ 
(since so is $\phi$, $\phi=G^F$ on $\pi^{-1}(\partial U)$, and $\phi\le G^F$)
and plurisubharmonic on $\bC^{k+1}\setminus(\pi^{-1}(\partial U)\cup\{0\})$,
and satisfies the mean value inequality at each point 
$Z_0\in\pi^{-1}(\partial U)\cup\{0\}$
on each complex line passing through $Z_0$ 
(since so does $\phi$, $\phi=G^F$ on $\pi^{-1}(\partial U)$, and $\phi\le G^F$). 
The function $\psi$ is also log-homogeneous (of order $1$).
By the log-homogeneity of both $G^F$ and $\psi$,
the function $\psi-G^F$ also descends to a function on $\bP^k$. 

Following the manner in \cite[\S 2.1]{BermanBoucksom10},
let us also denote by $\rd\rd^c\psi$ (resp.\ $\rd\rd^c G^F$)
the current {\em on} $\bP^k$ whose pullback under 
$\pi:\bC^{k+1}\setminus\{0\}\to\bP^k$
coincides (the genuine) $\rd\rd^c\psi$ (resp.\ the genuine $\rd\rd^c G^F$) 
{\em on} $\bC^{k+1}\setminus\{0\}$.
Then both the $k$-th {\itshape Bedford--Taylor wedge products}
$(\rd\rd^c\psi)^{\wedge k}$ and $(\rd\rd^c G^F)^{\wedge k}$ 
{\itshape on} $\bP^k$ of $\rd\rd^c\psi$ and $\rd\rd^c G^F$,
respectively, exist and
are {\itshape probability} measures on $\bP^k$. From the definition of $G^F$,
the latter probability measure $(\rd\rd^c G^F)^{\wedge k}$ on $\bP^k$
is nothing but the equilibrium measure 
$\mu_f$ of $f$. We claim that
$(\rd\rd^c\psi)^{\wedge k}=\mu_f$ on $\bP^k$;
for, by the definition of $\psi$, we have
$(\rd\rd^c\psi)^{\wedge k}=(\rd\rd^c G^F)^{\wedge k}=\mu_f$ 
on $\bP^k\setminus\overline{U}$, and
by $\supp\mu_f\subset J(f)$ and the vanishing $\mu_f(\partial U)=0$
in \eqref{eq:preimage}, we also have
$(\rd\rd^c\psi)^{\wedge k}(\overline{U})
=1-(\rd\rd^c\psi)^{\wedge k}(\bP^k\setminus\overline{U})
=1-\mu_f(\bP^k\setminus\overline{U})
=\mu_f(\overline{U})
=0$.
Hence the claim holds.

Once this claim is at our disposal, we have
$\psi-G^F\ge 0$ (indeed $=0$)
on $(\rd\rd^c\psi)^{\wedge k}$($=\mu_f$)-almost everywhere $\bP^k$,
and then by a version of the Bedford--Taylor domination principle 
(\cite[Corollary 2.5]{BEGZ10}; for a summary on the properties of
plurisubharmonic {\itshape weights} on {\itshape big} line bundles
over complex compact manifolds, which applies to
our $\omega_{\FS}$-plurisubharmonic functions {\itshape on} $\bP^k$,
see \cite[\S 2]{BermanBoucksom10}),
we have $\psi\ge G^F$ on $\bC^k\setminus\{0\}$. In particular
$\phi=\psi\ge G^F$ on $\pi^{-1}(U)$, which contradicts 
$U\cap\{\phi-G^F<0\}\neq\emptyset$. \qed

\begin{acknowledgement}
 This research was partially supported by JSPS Grant-in-Aid 
 for Scientific Research (C), 15K04924.
\end{acknowledgement}

\def\cprime{$'$}

\end{document}